%Paper: funct-an/9411007
%From: jdenis@esi.ac.at (Juriev Denis)
%Date: Tue, 22 Nov 1994 18:21:42 +0100 (MET)

\input amstex
\magnification=1200
\font\cyr=wncyr10
\font\cyb=wncyb10
\font\cyi=wncyi10
\documentstyle{amsppt}
\NoRunningHeads
\NoBlackBoxes
\define\vz{\text{\bf !`}}
\define\vy{\text{\bf !}}
\define\Aa{\frak A}
\define\Ab{\Aa^{\vz}}
\define\Ac{\Aa^{\widehat{\vz}}}
\define\Ad{(\Aa')^{\widehat{\vz}}}
\define\GL{\operatorname{GL}}
\define\gla{\operatorname{\frak g\frak l}}
\define\End{\operatorname{End}}
\define\aff{\operatorname{\frak a\frak f\frak f}}
\define\Mat{\operatorname{Mat}}
\define\Hom{\operatorname{Hom}}
\define\Map{\operatorname{Map}}
\define\Ki{\operatorname{Kir}}
\define\id{\operatorname{id}}
\define\codim{\operatorname{codim}}
\define\Vect{\operatorname{Vect}}
\define\eh{\operatorname{\text{\cyi e1g}}}
\topmatter
\title\cyb Harakteristiki par operatorov, gibridy Li, skobki
Puassona i nelinei0naya geometricheskaya algebra.
\endtitle
\author\cyr D.~V.~Yurp1ev
\endauthor
\address\tenpoint\cyr\newline
Otdel matematiki, Nauchno-issledovatelp1skie0\newline
institut sistemnyh issledovanie0 RAN
\endaddress
\email\nofrills\tenpoint\it E-mail address: \rm juriev\@systud.msk.su
\endemail
\address\tenpoint\rm\newline
Erwin Schr\"odinger Institut\newline
f\"ur Mathematische Physik,\newline
Wien, \"Osterreich (Austria)
\endaddress
\email\nofrills\tenpoint\it E-mail address: \rm jdenis\@esi.ac.at
\endemail
\date funct-an/yymmxxx\enddate
\endtopmatter
\document\cyr

Kak pokazyvaet istoriya matematiki, poisk reshenie0 "dikih", no vnutrenne
krasivyh zadach imeet podchas
bolp1shuyu znachimostp1, chem tochnye otvety k korrektno sformulirovannym
"ruchnym"  problemam. Odnoe0 iz ta\-kih "dikih" zadach yavlyaet\-sya
klas\-si\-fi\-ka\-tsi\-ya par linee0nyh operatorov v fiksirovannom linee0nom
topologicheskom (dazhe konechnomernom) prostranstve. Nesmotrya na yavnuyu
trans\-tsen\-dent\-nostp1 e1toe0 zadachi (i vo mnogom v silu e1to\-go),
"otblesk ei0 siyayushchego shlee0fa" -- analiz razlichnyh chislovyh i
algebraicheskih harakteristik par operatorov -- predstavlyaet nesomnennye0
interes kak s abst\-rakt\-noe0, tak i s prikladnoe0 tochek zreniya. I esli
{\cyi obratnaya zadacha teorii predstavlenie0\/} zaklyuchaet\-sya v
vo\-ssta\-nov\-le\-nii algebraicheskogo obp2ekta po zadannomu
se\-mee0\-s\-t\-vu
"predstavlyayushchih" ope\-ra\-to\-rov (matrits)(sm.napr. [1]), to ra\-zum\-no
poisk
razlichnyh struktur, yavlyayushchihsya algebraicheskimi, geometricheskimi i
geo\-met\-ro\-al\-geb\-rai\-ches\-ki\-mi harakteristikami
za\-dan\-no\-go se\-mee0\-s\-t\-va operatorov (matrits) otnositp1 k
{\cyi ras\-shi\-ren\-no\-mu tolkovaniyu\/} ob\-rat\-noe0 zadachi. Dannaya
rabota
posvyashchena pervonachalp1nomu izu\-che\-niyu ryada po\-dob\-nyh
ha\-rak\-te\-ris\-tik dlya pary operatorov.

Obshchie0 plan raboty takov. Para operatorov harakterizuet\-sya nekotorym
algebraicheskim obp2ektom, svoim "invariantom", gibridom Li (linee0nym
prostranstvom, osnashchi0nnym paroe0 soglasovannyh skobok Li). V opredeli0nnom
smysle ukazannye0 obp2ekt mozhno schitatp1 "stabilizatorom" infinitezimalp1nogo
dee0stviya drugogo algebraicheskogo obp2ekta, psevdogibrida Li (linee0nogo
prostranstva, osnashchi0nnogo paroe0 soglasovannyh struktur psevdoalgebr Li
[2])
na parah operatorov. Soglasovannostp1 struktur psevdoalgebr Li oznachaet, chto
lyubaya linee0naya kombinatsiya sootvet\-stvuyushchih strukturnyh funktsie0
snova zadai0t strukturu psevdoalgebry Li. Takim obrazom, na parah operatorov
zadano dvuparametricheskoe semee0stvo dee0stvie0 psevdoalgebr Li.
Bolee togo, na orbite psevdogibrida Li zadano dee0stvie nekotoroe0
bolp1shee0
psevdoalgebry Li (grifona dannoe0 orbity), poluchaemoe0 puti0m
dobavleniya k dee0stviyam psevdoalgebr Li, sootvet\-stvuyushchih
psevdogibridu Li, kommutativnogo semee0stva t.n. e1kvigibridnyh
variatsie0.

Takim
obrazom, para ope\-ra\-to\-rov mozhet bytp1 oharakterizovana kak strukturoe0
so\-ot\-vet\-s\-t\-vu\-yu\-shche\-go gibrida Li, tak i
geo\-met\-ri\-ches\-ki\-mi
svoe0\-s\-t\-va\-mi ei0 orbity i geometroalgebraicheskimi
svoe0\-s\-t\-va\-mi otvechayushchego ee0 grifona, a takzhe associirovannogo s
nim semee0stva izotopnyh lup [3,4]
(sm. takzhe [2]) i ih mulp1tialgebr Sabinina--Miheeva [5].

\definition{\cyb Opredelenie 1}

{\cyb A.} {\cyi Gibridom Li\/} \cyr nazyvaet\-sya linee0noe prostranstvo $V$,
snabzhi0nnoe pa\-roe0 soglasovannyh skobok Li $[\cdot,\cdot]'$ i
$[\cdot,\cdot]''$ (soglasovannostp1 oznachaet, chto lyubaya linee0naya
kombinatsiya e1tih skobok takzhe yavlyaet\-sya skobkoe0 Li) takih, chto
$$\aligned (\forall X,Y,Z\in V)\ &[[X,Z]',Y]''+[[X,Y]',Z]''+[[Z,Y]',X]''=\\
&[[X,Z]'',Y]'+[[X,Y]'',Z]'+[[Z,Y]'',X]'.\endaligned$$

{\cyb B {\cyr [6]}.} \cyr Para linee0nyh prostranstv $(V_1,V_2)$ nazyvaet\-sya
{\cyi
izotopicheskoe0 paroe0}, esli opredeleny otobrazheniya $m_1:V_2\otimes\wedge^2
V_1\mapsto V_1$ i $m_2:V_1\otimes\wedge^2V_2\mapsto V_2$ takie, chto skobki
$(X,Y)\mapsto[X,Y]_A=m_1(A,X,Y)$ ($X,Y\in V_1$, $A\in V_2$) i $(A,B)\mapsto
[A,B]_X=m_2(X,A,B)$ ($A,B\in V_2$, $X\in V_1$) yavlyayut\-sya soglasovannymi
skobkami Li pri razlichnyh znacheniyah podstrochnyh parametrov, bolee togo
$$\align
[X,Y]_{[A,B]_Z}=&\tfrac12([[X,Z]_A,Y]_B+[[X,Y]_A,Z]_B+[[Z,Y]_A,X]_B-\\
&[[X,Z]_B,Y]_A-[[X,Y]_B,Z]_A-[[Z,Y]_B,X]_A)\endalign$$ i $$\align
[A,B]_{[X,Y]_C}=&\tfrac12([[A,C]_X,B]_Y+[[A,B]_X,C]_Y+[[C,B]_X,A]_Y-\\
&[[A,C]_Y,B]_X-[[A,B]_Y,C]_X-[[C,B]_Y,A]_X),\endalign$$ gde ($X,Y,Z\in V_1$,
$A,B,C\in V_2$).
\enddefinition

\cyr Esli $\frak g$ --- algebra Li, to $(\frak g,\Bbb R)$ yavlyaet\-sya
izotopicheskoe0 paroe0.

\cyr Pustp1 $H$ -- linee0noe topologicheskoe prostranstvo, togda $(V_1,V_2)$
($V_1\simeq V_2\simeq\End(H)$) nadelyaet\-sya estestvennoe0 strukturoe0
izotopicheskoe0 pary ({\cyi operatornoe0 izotopicheskoe0 pary\/}:
$[X,Y]_A=XAY-YAX$, $[A,B]_X=AXB-BXA$ ($\forall X,Y\in V_1$, $\forall A,B\in
V_2$)). Konechnomernuyu operatornuyu izotopicheskuyu paru budem nazyvatp1
takzhe
{\cyi matrichnoe0 izotopicheskoe0 paroe0}.

\remark{\cyi Zamechanie 1} \cyr Kazhdomu gibridu Li $V$ otvechaet
izotopicheskaya para $(V,\Bbb R^2)$ s trivialp1nymi skobkami v $\Bbb R^2$ i
{\rm vice versa}.
\endremark

\cyr Pustp1 $(V_1,V_2)$ --- izotopicheskaya para, $\Aa$ -- linee0noe
podprostranstvo v $V_2$. Polozhim $\Ab=\{X\in V_1:\ (\forall
A,B\in\Aa)\ [A,B]_X\in\Aa\}$, $\Ac=\{X\in V_1:\ (\forall A,B\in\Aa)\
[A,B]_X=0\}$.

\proclaim{\cyb Teorema 1A} \cyi Dlya proizvolp1noe0 izotopicheskoe0 pary
$(V_1,V_2)$ i linee0nogo podprostranstva $\Aa\subseteq V_2$ pary $(\Ab,\Aa)$,
$(\Ac,\Aa)$ i $(\Ab/\Ac,\Aa)$ obladayut estestvennymi
strukturami izotopicheskih par; pri e1tom $(\Ab,\Aa)$ i $(\Ac,\Aa)$ dopuskayut
estestvennye vlozheniya v $(V_1,V_2)$.
\endproclaim

\cyr Pustp1 $(V_1,V_2)$ --- proizvolp1naya izotopicheskaya para i
$\Aa=\left<A,B\right>$ -- li\-nee0\-noe prostranstvo, natyanutoe na dva
e1lementa
$A$ i $B$ iz $V_2$. Polozhim $a(A,B)=\dim\Ac$, $a_0(A,B)=\dim \Ab/\Ac$.

\proclaim{\cyb Teorema 1B} \cyi $a_0(A,B)=0, 1, 2$. Esli $a_0(A,B)=1$, to
izotopicheskaya pa-\linebreak ra $(\Aa,\Ab/\Ac)$ yavlyaet\-sya estestvennoe0
izotopicheskoe0 paroe0, assotsiirovan-\linebreak noe0 s algebroe0 Li $\aff(\Bbb
R^1)$
gruppy affinnyh preobrazovanie0 pryamoe0. Esli\linebreak $a_0(A,B)=2$, to
izotopicheskaya para $(\Aa,\Ab/\Ac)$ sovpadaet s paroe0 Okubo $(\Bbb R^2,\Bbb
R^2)$ [1].
\endproclaim

\proclaim{\cyb Teorema 1V} \cyi Esli $(V_1,V_2)$ --- proizvolp1naya
izotopicheskaya para, $A$ i $B$ --- dva e1lementa iz $V_2$ i
$\Aa=\left<A,B\right>$, to $\Ac$ nadelyaet\-sya estestvennoe0 strukturoe0
gibrida Li, a imenno $(\forall X,Y\in\Ac)\ [X,Y]'=[X,Y]_A, [X,Y]''=[X,Y]_B$.
\endproclaim

\cyr Rassmotrim konechnomernye0 sluchae0.
\proclaim{\cyb Primer 1} \cyi Esli $A$ i $B$ dve matritsy $2\times 2$, to
$a(A,B)$ ravno 4, esli e1ti matritsy proportsionalp1ny, i 2 inache. Esli
$A=\left(\matrix a & b \\ c & d\endmatrix\right)$ i $B=\left(\matrix e & f \\
g & h\endmatrix\right)$ neproportsionalp1ny, to gibrid Li $\Ac$ porozhdaet\-sya
matritsami $$\left(\matrix de-ah & af-be \\ ag-ce & 0\endmatrix\right)\quad
\text{\cyi i}\quad\left(\matrix dg-ch & cf-bg \\ 0 & ag-ce\endmatrix\right)$$
i trivialen (t.e. obe skobki tozhdestvenno ravny nulyu). Esli $A$ i $B$
proportsionalp1ny, to gibrid Li $\Ac$ realizuet\-sya v prostranstve vseh
matrits $2\times 2$ i netrivialen (obe skobki proportsionalp1ny drug drugu i
imeyut vid $(X,Y)\mapsto XAY-YAX \text{ \cyi ili } XBY-YBX$
($X,Y\in\Mat_n(\Bbb R)$).
\endproclaim

\cyr Analogichno mozhet bytp1 obs\-chitan menee trivialp1nye0 sluchae0 $n=3$.
Zadacha nahozhdeniya $\Ab$ i $\Ac$ pri proizvolp1nom $n$ yavlyaet\-sya
linee0noe0 i legko programmiruemoe0.

\remark{\cyi Zamechanie 2} \cyr Esli $\Aa=\left<A,B\right>$ i
$\Aa'=\left<CAD,CBD\right>$, gde $C,D\in\GL(n,\Bbb R)$, to $\Ac$ i $\Ad$
izomorfny kak gibridy Li.
\endremark

\remark{\cyi Zamechanie 3} \cyr Pustp1 $F\in\Mat_n(\Bbb R)$, $\frak
h(F)=\left<F\right>^{\vy}=\{X\in\Mat_n(\Bbb R): [X,F]=0\}$, togda $\frak h(F)$
nadelyaet\-sya estestvennoe0 strukturoe0 gibrida Li so skobkami
$[X,Y]'=[X,Y]=XY-YX$ i $[X,Y]''=F[X,Y]=[X,Y]F=XFY-YFX$. Pustp1 $A\in\GL(n,\Bbb
R)$, $B\in\Mat_n(\Bbb R)$, $\Aa=\left<A,B\right>$, togda otobrazheniya
$X\mapsto XA$, $X\mapsto AX$ zadayut izomorfizm $\Ac$ na $\frak h(A^{-1}B)$ i
$\frak h(BA^{-1})$, so\-ot\-vet\-s\-t\-ven\-no. Pustp1 $A,B\in\Mat_n(\Bbb R)$ i
sushchestvuet matritsa $F$ takaya, chto $AF=B$ ($FA=B$),
$\Aa=\left<A,B\right>$,
togda otobrazhenie $X\mapsto XA$ ($X\mapsto AX$) zadai0t e1pimorfizm $\Ac$ na
$\frak h(F)$.
\endremark

\cyr Rassmotrim beskonechnomernye0 sluchae0. Pustp1 bazisnoe linee0noe
prostranstvo
--- prostranstvo mnogochlenov odnoe0 peremennoe0 $x$.
\proclaim{\cyb Primer 2} \cyi Esli $A=\frac{\partial}{\partial x}$ i $B=x$,
to gibrid Li porozhdaet\-sya dif\-fe\-ren\-tsi\-alp1\-ny\-mi operatorami $P$,
levye0 simvol kotoryh $P(x,\xi)$ udovletvoryaet differentsialp1nomu uravneniyu
$$P''_{x\xi}+xP'_x+\xi P'_{\xi}+P=0.$$
Zafiksiruem bazis $e_n=x^nf_n(x\frac{\partial}{\partial x})$ v prostranstve
reshenie0. Funkcii $f_n(\xi)$ imeyut vid $f(2k)=\frac{(2k-1)!!}{(2k+n)!!}$,
$f_n(2k+1)=0$. V bazise $e_n$ kommutatsionnye sootnosheniya v gibride Li
zadayut\-sya formulami
$$[e_m,e_n]^{\cdot}=\cases 0,& \text{ \cyr esli $m+n\in2\Bbb Z$,}\\
e_{n+m\pm1},& \text{ \cyr esli $m\in2\Bbb Z+1$, $n\in2\Bbb Z$,}\\
-e_{n+m\pm1},& \text{ \cyr esli $m\in2\Bbb Z$, $n\in2\Bbb Z+1$,}.\endcases$$
Znak plyus otvechaet skobke $[\cdot,\cdot]'$, a znak minus --- skobke
$[\cdot,\cdot]''$.
\endproclaim

\cyr Otmetim sleduyushchie0 fakt. Esli $V$ --- gibrid Li so skobkami
$[\cdot,\cdot]'$ i $[\cdot,\cdot]''$, to linee0noe prostranstvo $\frak k(V)=V
\oplus V$ nadelyaet\-sya strukturoe0 algebry Li s kommutatorom
$$\aligned&[(X_1,Y_1),(X_2,Y_2)]=\\
&\quad([X_1,X_2]''+\tfrac12([X_1,Y_2]'-[X_2,Y_1]'),
[Y_1,Y_2]'+\tfrac12([Y_1,X_2]''-[Y_2,X_1]'')).\endaligned$$

\definition{\cyb Opredelenie 2} \cyr Izotopicheskaya para $(V_1,V_2)$
nazyvaet\-sya {\cyi kon\-t\-ra\-gre\-di\-en\-t\-noe0}, esli $V_1=V_2^*$,
$V_2=V_1^*$ i dlya lyubyh $A,B\in V_2$, $X,Y\in V_1$ vypolnyaet\-sya tsepochka
ravenstv:
$$\left<[X,Y]_B,A\right>=-\left<[X,Y]_A,B\right>=\left<[A,B]_X,Y\right>=
-\left<[A,B]_Y,X\right>,$$ gde $\left<\cdot,\cdot\right>$ --- estestvennoe
sparivanie
$V_1$ i $V_2$.
\enddefinition

\remark{\cyi Zamechanie 4} \cyr Konechnomernaya operatornaya izotopicheskaya
para yavlyaet\-sya kontragredientnoe0.
\endremark

\cyr Opredelim dlya kontragredientnoe0 izotopicheskoe0 pary $(V_1,V_2)$
4--for\-mu $\Omega:\bigwedge^2V_1\otimes\bigwedge^2V_2\mapsto\Bbb R$:
%% FOLLOWING LINE CANNOT BE BROKEN BEFORE 80 CHAR
$$\Omega(A,B;X,Y)=\left<[X,Y]_B,A\right>=-\left<[X,Y]_A,B\right>=\left<[A,B]_X,Y\right>=
-\left<[A,B]_Y,X\right>.$$

\remark{\cyi Zamechanie 5} \cyr Vvidu nevyrozhdennosti sparivaniya $V_1\otimes
V_2\mapsto\Bbb R$ opredeleny operatory $R_{\Omega}^{(i)}:\bigwedge^2V_i\mapsto
\bigwedge^2V_i$ takie, chto
$$\left<R_{\Omega}^{(1)}(X\wedge Y),A\wedge B\right>=\Omega(A,B;X,Y)=
\left<X\wedge Y, R^{(2)}_{\Omega}(A\wedge B)\right>.$$
Pri e1tom operatory $R_{\Omega}^{(i)}$ sopryazheny drug drugu i
antiinvolyutivny:\linebreak $R^{(2)}_{\Omega}=(R^{(1)}_{\Omega})^*$,
$R^{(1)}_{\Omega}=(R^{(2)}_{\Omega})^*$, $(R^{(i)}_{\Omega})^2=-\id$.
\endremark

\remark{\cyi Zamechanie 6} \cyr Pustp1 $(V_1,V_2)$ --- proizvolp1naya
kontragredientnaya izotopicheskaya para, $A,B\in V_2$, $\Aa=\left<A,B\right>$,
togda $\Ac=\ker\Omega_{A,B}$, gde
$\Omega_{A,B}(\cdot,\cdot)=\Omega(A,B;\cdot,\cdot)$. Kak sledstvie, $\codim_
{\End(V_1)}\Ac$ --- chi0tno.
\endremark

\remark{\cyi Zamechanie 7} \cyr Pustp1 $A\in\GL(n,\Bbb R)$, $B\in\Mat_n(\Bbb
R)$, $\pi:X\mapsto XA \text{ \cyr ili } AX$ ($X\in\Mat_n(\Bbb R)$), togda
$\left.\Omega_{A,B}\right|_{\Mat_n(\Bbb R)/\Ac}=\pi_*\omega_{\Ki}(F)$ ---
obratnye0 obraz formy Kirillova [7,8] k orbite (ko)prisoedinennogo
predstavleniya gruppy Li $\GL(n,\Bbb R)$ tochki $F$ ($F=A^{-1}B \text{ \cyr
ili } BA^{-1}$) v e1toe0 tochke.
\endremark

\definition{\cyb Opredelenie 3}

\cyb A \cyr [2]. {\cyi Psevdoalgebroe0 Li\/} vektornyh polee0 na mnogoobrazii
$M$ na\-zy\-va\-et\-sya troe0ka $(\frak g, m, \tau)$, gde $\frak g$ ---
linee0noe
prostranstvo, $\tau\in\Hom(\frak g,\Vect(M))$,
$m\in\Map(M,\Hom(\bigwedge^2\frak
g,\frak g))$, prichi0m $m$ i $\tau$ soglasovany mezhdu soboe0, a imen\-no,
esli polozhitp1 $[X,Y]_a=m(a)(X,Y)$ ($X,Y\in\frak g$, $a\in M$), to
$[\tau(X),\tau(Y)](a)=\tau([X,Y]_a)$.

\cyb B \cyr. {\cyi Psevdogibridom Li\/} vektornyh polee0 na mnogoobrazii $M$
na\-zy\-va\-et\-sya linee0noe prostranstvo $\frak g$, osnashchi0nnoe dvumya
soglasovannymi strukturami psevdoalgebr Li $(\frak g, m', \tau')$ i
$(\frak g, m'', \tau'')$. Soglasovannostp1 oznachaet, chto
$[\tau'(X),\tau''(Y)](a)+[\tau''(X),\tau'(Y)](a)=\tau'([X,Y]''_a)+\tau''
([X,Y]'_a)$ ($X,Y\in\frak g$, $a\in M$), gde $[X,Y]'_a=m'(a)(X,Y)$ i
$[X,Y]''_a=m''(a)(X,Y)$.
\enddefinition

Esli $(\frak g, m, \tau)$ -- psevdoalgebra Li vektornyh polee0 na mnogoobrazii
$M$, to budem govoritp1, chto $\frak g$ dee0stvuet na $M$; analogichno, esli
$(\frak g, m', m'', \tau', \tau'')$ -- psevdogibrid Li vektornyh polee0 na
mnogoobrazii $M$, to budem govoritp1, chto $\frak g$ dee0stvuet na $M$.
Otmetim, chto esli $(\frak g, m', m'', \tau', \tau'')$ --- psevdogibrid Li, to
dlya lyubyh $\lambda$ i $\mu$ troe0ka $(\frak g, m_{\lambda,\mu},
\tau_{\lambda,
\mu})$ ($m_{\lambda,\mu}=\lambda m'+\mu m''$, $\tau_{\lambda,\mu}=\lambda\tau'+
\mu\tau''$) yavlyaet\-sya psevdoalgebroe0 Li; takim obrazom, dee0stvie
psevdogibrida Li na fiksirovannom mnogoobrazii zadai0t dvuparametricheskoe
semee0stvo psevdoalgebr Li $(\frak g, m_{\lambda,\mu}, \tau_{\lambda,\mu})$
vektornyh polee0 na ni0m.

\cyr Pustp1 $(V_1,V_2)$ --- proizvolp1naya izotopicheskaya para. Polozhim
$W_i=V_i\oplus V_i$.

\proclaim{\cyb Teorema 2A}
\cyi $V_1$ dee0s\-t\-vu\-et kak psevdogibrid Li na
$W_2$ (i naoborot, $V_2$ dee0s\-t\-vu\-et kak psevdogibrid Li na $W_1$):
$$\aligned\tau'(X)(A,B)=([A,B]_X,0),&\quad\tau''(X)(A,B)=(0,[A,B]_X)\\
\tau'(A)(X,Y)=([X,Y]_A,0),&\quad\tau''(A)(X,Y)=(0,[X,Y]_A)\endaligned$$
($X,Y\in\frak f_1$, $A,B\in\frak f_2$); pri e1tom
$$\aligned
[X,Y]'_{(A,B)}=[X,Y]_B,&\quad[X,Y]''_{(A,B)}=[X,Y]_A\\
[A,B]'_{(X,Y)}=[A,B]_Y,&\quad[A,B]''_{(X,Y)}=[A,B]_X
\endaligned$$
\endproclaim

\remark{\cyi Zamechanie 8} \cyr Esli konechnomernye0 psevdogibrid Li $\frak g$
dee0stvuet na konechnomernom mnogoobrazii $M$, to kazhdaya psevdoalgebra Li
vektornyh polee0 na $M$ iz sootvet\-stvuyushchego dvuparametricheskogo
semee0stva psevdoalgebr Li e1ks\-po\-nen\-tsi\-ru\-et\-sya do
psev\-do\-(kva\-zi)\-grup\-py Li
preorazovanie0 $M$ [2-4,9].
\endremark

\definition{\cyb Opredelenie 4}

\cyb A. \cyr {\cyi Orbitoe0\/} psevdoalgebry Li $(\frak g, m, \tau)$ vektornyh
polee0 na mnogoobrazii $M$ nazyvaet\-sya minimalp1noe podmnogoobrazie $N$
takoe, chto $\tau(\frak g)\subseteq\Vect(N)$.

\cyb B. \cyr {\cyi Orbitoe0\/} psevdogibrida Li $(\frak g, m', m'', \tau',
\tau'')$ vektornyh polee0 na mnogoobrazii $M$ nazyvaet\-sya minimalp1noe
podmnogoobrazie $N$ takoe, chto\linebreak
$\tau'(\frak g)+\tau''(\frak g)\subseteq
\Vect(N)$.
\enddefinition

\remark{\cyi Zamechanie 9} \cyr Dlya proizvolp1noe0 izotopicheskoe0 pary
$(V_1,V_2)$ prostranstvo $W_2$ dopuskaet $(V_1, m_{\lambda,\mu}, \tau_
{\lambda,\mu})$--e1kvivariantnoe rassloenie $\Cal P_{\lambda,\mu}$ nad
bazoe0 $V_2$ s proektsiee0
$p:W_2\mapsto V_2$, gde $p(A,B)=\lambda A+\mu B$, sloe0 $U(C)$ nad tochkoe0
$C$ izomorfen $V_2$ dlya lyubogo $C$. Na lyubom sloe $U(C)$ dee0stvie
psevdoalgebry Li $(V_1, m_{\lambda,\mu}, \tau_{\lambda,\mu})$ redutsiruet\-sya
do dee0stviya algebry Li $\frak l(C)$, realizuemoe0 v prostranstve $V_1$ s
pomoshchp1yu $[\cdot,\cdot]_C$, $U(C)\simeq\frak l^*(C)$. Takim obrazom,
orbity dee0stviya psevdoalgebry Li $(V_1, m_{\lambda,\mu}, \tau_{\lambda,\mu})$
na $W_2$ sovpadayut s orbitami koprisoedini0nnyh predstavlenie0 algebr Li
$\frak l(C)$.
\endremark

\remark{\cyi Zamechanie 10}
\cyr Dlya matrichnoe0 izotopicheskoe0 pary $(V_1,V_2)$ orbity
dee0\-s\-t\-viya
psevdogibrida Li $V_1$ na $W_2$ sovpadayut s orbitami dee0stviya $\GL(n,\Bbb R)
\times\GL(n,\Bbb R)$, opisannogo v zamechanii 2.
\endremark

\cyr Takim obrazom, kasatelp1noe prostranstvo k orbite psevdogibrida Li $V_1$
na $W_2$ v proizvolp1noe0 tochke dlya matrichnoe0 izotopicheskoe0 pary
ne sovpadaet s pryamoe0 summoe0 kasatelp1nyh prostranstv k orbitam
so\-ot\-vet\-
s\-t\-vu\-yu\-shchih psevdoalgebr Li v ukazannoe0 tochke.

\definition{\cyb Opredelenie 5} \cyr Pustp1 $(V_1,V_2)$ --- proizvolp1naya
izotopicheskaya para, $\Cal O$ --- orbita dee0stviya psevdogibrida Li
$V_1$ na $W_2$. Opredelim {\cyi e1kvigibridnoe sloenie\/} $\Cal F_{\eh}$
na orbite $\Cal O$ sleduyushchim obrazom: sloi $\Cal F_{\eh}$ sostoyat iz
par $(A,B)$ s sovpadayushchimi gibridami Li $\Ac$ ($\Ac=\left<A,B\right>$).
\enddefinition

\remark{\cyi Gipoteza} \cyr Dlya proizvolp1noe0 izotopicheskoe0 pary
$(V_1,V_2)$ kasatelp1noe pro\-s\-t\-ran\-s\-t\-vo k orbite dee0stviya
psevdogibrida Li
$V_1$ na $W_2$ v proizvolp1noe0 tochke yavlyaet\-sya summoe0 kasatelp1nyh
prostranstv k orbitam so\-ot\-vet\-s\-t\-vu\-yu\-shchih psevdoalgebr Li i k
e1kvigibridnomu
sloeniyu v ukazannoe0 tochke.
\endremark

\cyr Dannaya gipoteza spravedliva dlya matrichnoe izotopicheskoe0 pary. Bo\-lee
togo, v e1tom sluchae e1kvigibridnoe sloenie porozhdaet\-sya kommutativnym
semee0stvom variatsie0, kotorye my budem nazyvatp1 {\cyi e1kvigibridnymi
variatsiyami}. A imenno, dlya proizvolp1nogo $X\in\Ac$ ($\Ac=\left<A,B\right>$)
opredelim variatsii $\delta_{\eh}(X)$ sleduyushchim obrazom: $\delta_{\eh}(X)
(A,B)=(AXA+AXB, AXB+BXB)$. Netrudno pokazatp1, v e1tom sluchae v
proizvolp1noe0 tochke
orbity psevdogibrida Li imeet mesto razlozhenie kasatelp1nogo prostranstva k
orbite v {\cyi pryamuyu summu\/} kasatelp1nyh prostranstv k orbitam
sootvet\-stvuyushchih psevdoalgebr Li i e1kvigibridnyh variatsie0
(kotorye obrazuyut kommutativnuyu algebru Li). Kak sledstvie, na
orbite psevdogibrida Li $V_1$ zadano dee0stvie psevdoalgebry Li
$\widehat W_1=W_1\oplus\Ac$, gde dee0stvie $W_1$ opredelyaet\-sya dee0stviyami
psevdoalgebr Li $(V_1, m', \tau')$ i $(V_2, m'', \tau'')$, a
dee0stvie $\Ac$ zadai0t\-sya e1kvigibridnymi variatsiyami. Ukazannuyu
psevdoalgebru Li $\widehat W_1$ budem nazyvatp1 {\cyi grifonom\/}
sootvet\-stvuyushchee0 orbity psevdogibrida Li.

\remark{\cyi Zamechanie 11} \cyr V konechnomernom sluchae dee0stvie
grifona na orbite e1ks\-po\-nen\-tsi\-ru\-et\-sya do dee0stviya
so\-ot\-vet\-s\-t\-vu\-yu\-shchee0 psevdo(kvazi)gruppy Li [2-4,9].
\endremark

\remark{\cyi Zamechanie 12} \cyr Razlozhenie kasatelp1nogo
prostranstva k orbite $\Cal O$ psevdogibrida Li $V_1$ v pryamuyu
summu kasatelp1nyh prostranstv k orbitam
psevdoalgebr Li $(V_1, m_{\lambda',\mu'}, \tau_{\lambda',\mu'})$ i
$(V_2, m_{\lambda'',\mu''}, \tau_{\lambda'',\mu''})$ ($\left|\matrix
\lambda' & \lambda'' \\ \mu' & \mu''\endmatrix\right|\ne 0$) i
e1kvigibridnyh variatsie0 oznachaet, chto e1kvigibridnoe sloenie i
rassloeniya $\Cal P_{\lambda',\mu'}$, $\Cal P_{\lambda'',\mu''}$
zadayut {\cyi setp1\/} [10] na orbite $\Cal O$ psevdogibrida Li
$V_1$. Otmetim takzhe, chto $(n+1)$-ka
$(\Cal P_{\lambda_1,\mu_1},\ldots\Cal P_{\lambda_n,\mu_n},\Cal F_{\eh})$
($\left|\matrix \lambda_i & \lambda_j \\ \mu_i &
\mu_j\endmatrix\right|\ne 0$ esli $i\ne j$) zadai0t nekotorye0
obp2ekt na orbite $\Cal O$, svoe0stva kotorogo s\-hodny so svoe0stvami
{\cyi mnogomernoe0 $n$--tkani\/} [11,12].
\endremark

\remark{\cyi Zamechanie 13} \cyr Lupa, opredelyaemaya grifonom
orbity psevdogibrida Li matrichnoe0 izotopicheskoe0 pary (ili
sootvet\-stvuyushchee0 psev\-do\-(kva\-zi)\-grup\-poe0 Li), i tochkoe0
orbity $(A,B)$ [3,4] dopuskaet vlozhenie v sebya gruppy
$\frak L_{\lambda,\mu}(C)$,
otvechayushchee0
algebre Li $\frak l(C)$ $(C=\lambda A+\mu B)$, dlya lyubyh $\lambda$
i $\mu$, a takzhe kommutativnoe0 gruppy,
otvechayushchee0 e1kvigibridnym variatsiyam.
\endremark

\cyr Opredelim rassloenie $\Cal P_{\widehat{\vz}}$ nad orbitoe0
$\Cal O$ psevdogibrida Li, sloe0 kotorogo nad tochkoe0 $(A,B)$
sovpadaet s $\left<A,B\right>^{\widehat{\vz}}$. V rassloenii $\Cal P_
{\widehat{\vz}}$ zadana svyaznostp1 $\nabla$ takaya, chto
$$\aligned&\nabla_{\tau'(Z_1)+\tau''(Z_2)+\delta_{\eh}(Z_0)}((A,B),X)=\\
&\qquad((\tau'(Z_1)+\tau''(Z_2)+\delta_{\eh}(Z_0))(A,B),
[Z_1,X]_B+[Z_2,X]_A),\endaligned$$
gde $Z_1, Z_2\in V_2$, $Z_0, X\in\left<A,B\right>^{\widehat{\vz}}$.
Rassmotrim takzhe dvoe0stvennoe rassloenie $\Cal P^*_{\widehat{\vz}}$
nad $\Cal O$ so sloem $({\left<A,B\right>^{\widehat{\vz}}})^*$ nad
tochkoe0 $(A,B)$, v kotorom opredelena dvoe0stvennaya svyaznostp1
$\nabla^*$.
Kovariantnye proizvodnye v $\Cal P^*_{\widehat{\vz}}$ vdolp1
e1lementov grifona orbity za\-my\-ka\-yut\-sya do dee0stviya
nekotoroe0 psevdoalgebry Li, kotoruyu budem nazyvatp1
{\cyi kanonicheskim rasshireniem grifona.}

\remark {\cyi Zamechanie 14} \cyr Dlya matrichnoe0 izotopicheskoe0
pary svyaznostp1 $\nabla^*$ v rassloenii $\Cal P^*_{\widehat{\vz}}$
zadai0t podnyatie dee0stviya dublya $\GL(n,\Bbb R)\times\GL(n,\Bbb R)$
na baze do ego svobodnogo i tranzitivnogo dee0stviya
$\gla(n,\Bbb R)\oplus\gla(n,\Bbb R)$ na totalp1nom prostranstve.
Kak sledstvie, na totalp1nom prostranstve rassloeniya $\Cal P^*_
{\widehat{\vz}}$ zadana puassonova struktura dublya Gee0zenberga
[13,14].
\endremark

\definition{\cyb Opredelenie 6 \cyr (sr. [2])}
\cyr Psevdoalgebra Li $(\frak g, m, \tau)$ nazyvaet\-sya
{\cyi ga\-milp1\-to\-no\-voe0}, es\-li so\-ot\-vet\-s\-t\-vu\-yu\-shchee
mnogoobrazie $M$ --
puassonovo i zadano linee0noe otobrazhenie $\hat\tau\in\Hom(\frak g,\Cal O(M))$
takoe, chto $\tau(X)F=\{\hat\tau(X),F\}$ ($X\in\frak g$, $F\in\Cal O(M)$.
Esli $(\frak g, m, \tau)$ -- gamilp1tonova psevdoalgebra Li vektornyh polee0
na puassonovom mnogoobrazii $M$, to budem govoritp1, chto $\frak g$ dee0stvuet
na $M$ gamilp1tonovo.
\enddefinition

\proclaim{\cyb Teorema 2B} \cyr Puassonova struktura na
$\Cal P^*_{\widehat{\vz}}$ invariantna otnositelp1no
dee0\-s\-t\-viya
kanonicheskogo rasshireniya grifona, yavlyayushchegosya
gamilp1tonovym.
\endproclaim

\cyr
Avtor nadeet\-sya, chto dannaya zametka posluzhit stimulom dlya
us\-ta\-nov\-le\-niya bolee tesnyh svyazee0 mezhdu klassicheskimi
napravleniyami funk\-tsio\-nalp1\-no\-go analiza s odnoe0 storony i
sovremennymi issledovaniyami po nelinee0nym skobkam Puassona
(A.Vae0nshtee0n--M.V.Karasi0v--V.P.Maslov [2]), nelinee0noe0 geometricheskoe0
algebre (L.V.Sa\-bi\-nin--P.O.Miheev [15,\linebreak 16], sm. takzhe
[17]) i
bes\-ko\-nech\-no\-mer\-noe0
simplekticheskoe0 geo\-met\-rii s drugoe0.

Avtor blagodarit Laboratoriyu Teoreticheskoe0 Fiziki Vysshee0 Normalp1noe0
Shkoly (Parizh) za isklyuchitelp1nuyu atmosferu, v kotoroe0 voz\-nik zamysel
e1toe0 raboty, i Mezhdunarodnye0 Institut Matematicheskoe0 Fiziki imeni
E1rvina Shri0dingera (Vena), gde rabota byla zavershena, za podderzhku i
lyubeznoe gostepriimstvo.

\

\

\head{\cyr SPISOK LITERATURY}\endhead

\

\roster
\item"[1]" {\rm Juriev D., Topics in hidden symmetries, E--print (LANL Archive
on Theor. High Energy Phys.): {\it hep-th/9405050} (1994).}
\item"[2]" {\cyr Karasi0v M.V., Maslov V.P., Ne\-li\-nee0\-nye
skob\-ki Pu\-as\-so\-na:
geo\-met\-riya i kvan\-to\-va\-nie, M., Na\-u\-ka, 1992.}
\item"[3]" {\cyr Miheev P.O., O lu\-pah
pre\-ob\-ra\-zo\-va\-nie0. V sb. "Ne\-ko\-to\-rye pri\-lo\-zhe\-niya
dif\-fe\-ren\-tsi\-alp1\-noe0 geo\-met\-rii". M., 1985, S.85-93. Ruk. dep. v
VINITI 4531-85Dep.}
\item"[4]" {\rm Mikheev P.O., Quasigroups of transformations, Trans. Inst.
Phys. Estonian Acad. Sci. 1990. V.66. P.54-66.}
\item"[5]" {\cyr Sabinin L.V., Miheev P.O., Ob
in\-fi\-ni\-te\-zi\-malp1\-noe0 te\-o\-rii
lo\-kalp1\-nyh ana\-li\-ti\-ches\-kih lup, Dok\-la\-dy AN SSSR, 1988, T.297,
S.801-805.}
\item"[6]" {\rm Juriev D., Topics in isotopic pairs and their representations,
E--print (Texas Archive on Math. Phys.): {\it mp\_arc/94-267} (1994).}
\item"[7]" {\cyr Kirillov A.A., E1le\-men\-ty te\-o\-rii
pred\-stav\-le\-nie0, M., Na\-u\-ka,
1971.}
\item"[8]" {\cyr Fomenko A.T., Simp\-lek\-ti\-ches\-kaya
geo\-met\-riya. Me\-to\-dy i
pri\-lo\-zhe\-niya, M., Izd-vo MGU, 1988.}
\item"[9]" {\rm Batalin I., Quasigroup construction and first class
constraints.
J. Math. Phys. 1981. V.22. P.1837-1850.}
\item"[10]" {\cyr Ba\-zy\-lev V.T., Kuzp1\-min M.K., Sto\-lya\-rov
A.V., Se\-ti na mno\-go\-ob\-ra\-zi\-yah. V sb. "Problemy geometrii
12". M., VINITI, 1981, S.97-125.}
\item"[11]" {\rm Gol'dberg V., Theory of multicodimensional
$(n+1)$--webs. Kluwer Acad. Publ., Dordrecht, 1988.}
\item"[12]" {\rm Akivis M.A., Shelekhov A.M., Geometry and algebra of
multidimensional 3-webs. Kluwer Acad. Publ., Dordrecht, 1992.}
\item"[13]" {\rm Semenov-Tian-Shansky M.A., Dressing transformations
and Poisson-Lie\linebreak group actions. Publ. RIMS Kyoto, 1985, V.21,
P.1237-1260.}
\item"[14]" {\rm Alekseev A.Yu., Malkin A.Z., Symplectic structures
associated to Lie-Pois\-son groups. Commun. Math. Phys. 1994.
V.***. P.***-***; E--print (LANL Archive on High Energy Phys):
{\it hep-th/9303038} (1993).}
\item"[15]" {\cyr Sabinin L.V., O ne\-li\-nee0\-noe0
geo\-met\-ri\-ches\-koe0 al\-geb\-re. V sb. "Tka\-ni
i kva\-zi\-grup\-py". Ka\-li\-nin (Tverp1), 1988, S.32-37.}
\item"[16]" {\cyr Sabinin L.V., Miheev P.O., Glad\-kie
kva\-zi\-grup\-py i geo\-met\-riya. V sb.
"Prob\-le\-my geo\-met\-rii 20". M., VINITI, 1988, S.75-110.}
\item"[17]" {\cyr Sabinin L.V., Me\-to\-dy ne\-as\-so\-tsia\-tiv\-noe0
al\-geb\-ry v
dif\-fe\-ren\-tsi\-alp1\-noe0 geo\-met\-rii, Pri\-lozh. k russk. per. kn.
Ko\-ba\-yasi S., No\-mi\-dzu
K., Os\-no\-vy dif\-fe\-ren\-tsi\-alp1\-noe0 geo\-met\-rii. T.1. M.,
Na\-u\-ka, 1982.}
\endroster
\enddocument